\documentclass[submission,copyright,creativecommons]{eptcs}
\providecommand{\event}{GASCOM 2024} % Name of the event you are submitting to

\usepackage{iftex}

\ifpdf
  \usepackage{underscore}         % Only needed if you use pdflatex.
  \usepackage[T1]{fontenc}        % Recommended with pdflatex
\else
  \usepackage{breakurl}           % Not needed if you use pdflatex only.
\fi

\usepackage{amsmath, amsfonts, amsthm, amssymb}
\usepackage{latexsym, mathrsfs} 
\usepackage{mathtools}
\usepackage{enumitem}
\usepackage{enumerate}
\usepackage{todonotes}

\usepackage{tikz}
\usetikzlibrary{patterns}
\usetikzlibrary{fadings}
\tikzfading[name=fade out, inner color=transparent!0, outer color=transparent!90]
\usepackage[labelfont=normalfont,labelformat=simple]{subcaption}
\usepackage{graphicx}
\usepackage[vcentermath]{youngtab}
\usepackage{ytableau}

\usepackage[capitalise]{cleveref}

\newtheorem{theorem}{Theorem}[section]

\title{A Symmetry Property of Christoffel Words}

\author{Yan Lanciault
\institute{LACIM, Université du Québec à Montréal,\\ Montréal, Québec}
\email{\href{mailto:lanciault.yan@courrier.uqam.ca}{lanciault.yan@courrier.uqam.ca}}
\and
Christophe Reutenauer\thanks{Christophe Reutenauer was partially supported by NSERC}
\institute{LACIM, Université du Québec à Montréal,\\ Montréal, Québec}
\email{\href{reutenauer.christophe@uqam.ca}{reutenauer.christophe@uqam.ca}}
}
\def\titlerunning{A Symmetry Property of Christoffel Words}
\def\authorrunning{Y. Lanciault, C. Reutenauer}

\begin{document}

\maketitle

%\tableofcontents

\begin{abstract} 
Motivated by the theory of trapezoidal words, whose sequences of cardinality of factors by length are symmetric, we introduce a bivariate variant of this symmetry. We show that this symmetry characterizes Christoffel words, and prove other related results. 
\end{abstract}

\section{Introduction}

{\em Trapezoidal words} were considered by Aldo de Luca in \cite{dL}; for such a word, $w$ say, of length $n$, the graph of the discrete function 
$\{0,1,\ldots,n\}\to \mathbb N$, giving the number of factors of length $k$ of $w$ is an isosceles trapezoid, with successive values 
$1,2,\ldots,J,J+1,\ldots,J+1,J,\ldots,2,1$.
%the definition is somewhat technical: with a word $w$ he associates two parameters $K$ and $L$, and then a word is trapezoidal if and only if $K+L$ is equal to the length of $w$; 
He showed that Sturmians words are trapezoidal, but the converse does not necessarily hold. The terminology ``trapezoidal" was introduced by Flavio d'Alessandro in 
\cite{DA}, who studied these words, giving in particular a condition for which a trapezoidal word is Sturmian. In \cite{BDLF}, Michelangelo Bucci, 
Alessandro De Luca and Gabriele Fici gave many equivalent conditions for a word to be trapezoidal; one of them is that the number of factors
of length $k$ is at most $k+1$ (also see the work of Florence Levé and Patrice Séébold \cite{LS}, and that of Mira-Cristiana Anisiu and Julien 
Cassaigne \cite{AC}). Remind that a factor of a word is a contiguous subword.

A remarkable property of trapezoidal words is, as mentioned above, that the sequence of the lengths of the factors of these words, from length 0 to 
length $n$, is symmetric. We may call such a word {\em factor-symmetric}.

In the present work, we present a generalization of this symmetry property. Let $w$ be a word over the alphabet $\{a,b\}$, with $p$ occurrences of the letter 
$a$ and $q$ occurrences of the letter $b$; in other words, the Parikh image of $w$ is $(p,q)$. We say that $w$ is {\em strongly factor-symmetric} if 
for any $i,j$, $w$ has as many distinct factors with Parikh image $(i,j)$ as distinct factors of Parikh image $(p-i,q-j)$. Note that in that case, the notion of symmetry does not necessarily mean invariant under reversal.

We show that each Christoffel word is strongly factor-symmetric (Theorem \ref{ECwSS}). Conversely, each finite primitive Sturmian word which is 
strongly 
factor-symmetric is a Christoffel word (Theorem \ref{SSSwC}). Note that $aabb$ is strongly factor-symmetric, so that the hypothesis ``Sturmian" is 
not superfluous.

These results are interesting, in part because one obtain a characterization of Christoffel words among all Sturmian words. Indeed, in the 
literature there exist many characterizations of conjugate of Christoffel words (\cite{Ch2, MRS, BR, R1, MRRRS, R2}), which do not distinguish between Christoffel words 
and their conjugates. However, another notable charaterization of Christoffel words is that a Sturmian word is a Christoffel word if and only if it is a Lyndon word 
\cite{BdL}, if and only if it is unbordered \cite{Ch1} (see also \cite{HN}).
%(plein de references ICI).

Concerning nonprimitive words, we show that if $w$ is
%that powers on Christoffel words strongly factor-symmetric (Theorem ref{}; whose proof uses other tehcniques), and that conversely, a  
a nontrivial power of a primitive word $u$, then $w$ is strongly factor-symmetric if and only if $u$ is a Christoffel word (Theorem \ref{SSiffuCw}). The hypothesis 
``Sturmian" is not necessary here. In particular, $(aabb)^2$ is not strongly factor-symmetric.

As a byproduct, we obtain that, with the notation of the previous paragraph, that $w$ is factor-symmetric if and only if $u$ is the conjugate of some Christoffel 
word (Theorem \ref{conjugate}).

Concerning the strong factor symmetry of a Christoffel word $w$, we give an explicit bijection between the factors of $w$ of Parikh image $(i,j)$ and those of 
Parikh image $(p-i,q-j)$ (Theorem \ref{bijection}); it relies on the notion of {\em attractor} and {\em circular attractor} \cite{MRRRS}. Moreover, the support of the function of pairs of integers that counts the numbers of factors of $w$ for each Parikh image, which is a subset of the discrete plane, is the set of integer points on the two paths defined by $w$ and its reversal $\tilde w$, and between them (Theorem \ref{support}): see, for example (\ref{aabab-array}) and Figure \ref{Cword2}.

This work was partially supported by NSERC, Canada. 

\section{Christoffel words and Sturmian words}

Among several equivalent definitions of Christoffel words, we choose the following: a Christoffel word on the alphabet $\{a,b\}$ is either $a$ or $b$, or a word of the form $amb$ or $bma$, such that $m$ is a palindrome, and $w$ is a product of two palindromes. For other characterizations, see for example the book of the second author \cite {R1}. Christoffel words are primitive, that is, are not equal to a nontrivial power of another word.

It is known that the factorization into two palindromes is unique, and it is called the {\em palindromic factorization}.

Given a word $w$, we define the function $\delta_w: \mathbb N^2\to \mathbb N$ by $\delta_w(i,j)=$ the number of factors of $w$ whose Parikh image is $(i,j)$. We say that a word $w$ of Parikh image $(p,q)$ is {\em strongly factor-symmetric} if for any $i,j$, $
\delta_w(i,j)=\delta_w(p-i,q-j)$. For example, the distinct factors of the Christoffel word $aabab$ are $1,a,b,aa,ab,ba,aab,aba,bab,aaba,abab,aabab$ so that $
\delta_w$ is represented by the array whose $i,j$-coordinate is $\delta_w(i,j)$ (coordinates are as in the Cartesian plane, and this array is embedded in the plane):
\begin{equation}\label{aabab-array}
%\left[
\begin{array}{cccccc}
0&1&1&1\\
1&2&2&1\\
1&1&1&0
\end{array}
%\right]
\end{equation}
This array has a central symmetry, which means that $w$ is strongly factor-symmetric. We call this array the {\em factor array} of $w$.

A word $w$ is called {\em factor-symmetric} if the sequence of length of factors, which turns out to be $\sum_{i+j=k}\delta_w(i,j)$, $k=0,\ldots,|w|$ is symmetric; in other words, $w$ has as many factors of length $i$ as factors of length $n-i$, for all $i$, with $n=|w|$.
Clearly, a strongly factor-symmetric word is factor-symmetric. 

Trapezoidal words are factor-symmetric words (\cite{dL} Proposition 4.7, \cite{BDLF} Definition 2.5); and conversely, each factor-symmetric word $w$ is trapezoidal: indeed, if $|w|=n$, then $w$ has $n-i+1$ occurrences of factors of length $i$, so that it has at most $n-i+1$ such factors; but the factor symmetry implies that it has at most $n-i+1$ factors of length $n-i$, and hence it is trapezoidal by the cited proposition.

\section{Main results}

\begin{theorem}\label{ECwSS} Each Christoffel word is strongly factor-symmetric.
\end{theorem}

We have a converse. Note that a Sturmian word is a factor of a Christoffel word.

\begin{theorem}\label{SSSwC} If the support of $\delta_w$ is symmetric (and in particular, if $w$ is strongly factor-symmetric) and if $w$ is primitive and Sturmian, then $w$ is a Christoffel word.
\end{theorem}

Note that the factor array of the word $aabb$ is
$$
\begin{array}{cccccc}
1&1&1\\
1&1&1\\
1&1&1
\end{array}
$$
which has a central symmetry, so that $aabb$ is strongly factor-symmetric; this word is not a Christoffel word, but is not Sturmian either, since $aa$ and $bb$ 
cannot be both factors of a Sturmian word.

\begin{theorem}\label{SSiffuCw} Let $w=u^k$, $u$ primitive, $k\geq 2$. Then $w$ is strongly factor-symmetric if and only if $u$ is a Christoffel word.
\end{theorem}

Here, the hypothesis ``factor-symmetric" suffices for the ``only if" part. And the hypothesis ``Sturmian" is no more necessary.

\section{Byproducts}
An {\em attractor } of a word $w=w_1\cdots w_n$, with $w_i$ letters of the alphabet, is a subset $K$ of $\{1,\cdots,n\}$ such that every factors of $w$ has an occurrence that meets one of the letters indexed by one of the numbers in $K$. A {\em circular attractor } is defined similarly, but with the notion of circular factors, that is factors of a conjugate of $w$.   
Using theses concepts, we have a bijection that explains Theorem \ref{ECwSS}.

\begin{theorem}\label{bijection} Let $w=uv$ be a Christoffel word of length $n$ with its palindromic factorization. Suppose $k, 0\leq k\leq n $. Consider all factors of length $k$ of $w$ 
that intersect the cut of the factorization, and order them from left to right: $f_1,f_2,\ldots,f_r$. Consider all factors of length $n-k$ of $w$ that intersect this cut, 
and order them from right to left: $g_1,g_2,\ldots,g_s$. Then $r=s$, the words $f_i$ are distinct, the words $g_i$ are distinct, and the mapping $f_i\mapsto g_i$ is 
a bijection from the set of factors of length $k$ of $w$ to 
the set of factors of length $n-k$ of $w$, which complements the Parikh image $\gamma(w)$ of $w$; that is: $\gamma(f_i)+\gamma(g_i)=\gamma(w)$.
\end{theorem}

An example: let $w=aababab$, $u\cdot v=aa\cdot babab$, $k=4$, $f_1=aaba,f_2=abab,f_3=baba$, $g_1=bab,g_2=aba,g_3=aab$.

In the following, with each word on the alphabet $\{a,b\}$, we associate the path in the discrete plane starting from the origin, where $a$ represents an horizontal step towards East, and $b$ a vertical step towards North.

\begin{theorem}\label{support} Let $w$ be a lower Christoffel word, $\tilde w$ the corresponding upper Christoffel word, and $S_w$ the set of integer points on the paths corresponding to $w$ and $\tilde w$. Then $S_w$ is the support of the function $\delta_w$.
\end{theorem}

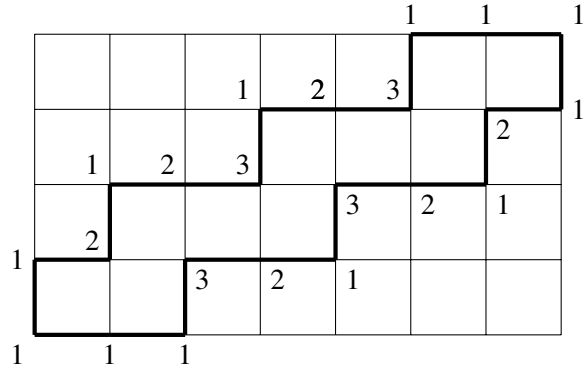
\begin{figure}
\centering
\begin{tikzpicture}
\draw (0,0) grid (7,4); 
%\draw[thick] (0,0) -- (7,4);

\draw [ultra thick] (0,0) -- (1,0); \draw  [ultra thick] (0,0) -- (0,1);
\draw [ultra thick] (7,3) -- (7,4); \draw  [ultra thick] (6,4) -- (7,4);

\draw [ultra thick] (1,0) -- (2,0) ; 
\draw [ultra thick] (2,0) -- (2,1) ;
\draw [ultra thick] (2,1) -- (4,1);
\draw [ultra thick] (4,1) -- (4,2);
\draw [ultra thick] (4,2) -- (6,2);
\draw [ultra thick] (6,2) -- (6,3);
\draw [ultra thick] (6,3) -- (7,3);

\draw [ultra thick]  (0,1) -- (1,1); 
\draw [ultra thick] (1,1) -- (1,2);
\draw [ultra thick] (1,2) -- (3,2);
\draw [ultra thick] (3,2) -- (3,3);
\draw [ultra thick] (3,3) -- (5,3);
\draw [ultra thick] (5,3) -- (5,4);
\draw [ultra thick]  (5,4) -- (6,4);

\draw (0,0) node[below left] {$1$};
\draw (1,0) node[below] {$1$};
\draw (2,0) node[below] {$1$};
\draw (0,1) node[left] {$1$};
\draw (1,1) node[above left] {$2$};
\draw (1,2) node[above left] {$1$};
\draw (2,1) node[below right] {$3$};
\draw (3,1) node[below right] {$2$};
\draw (2,2) node[above left] {$2$};
\draw (4,1) node[below right] {$1$};
\draw (3,2) node[above left] {$3$};
\draw (4,2) node[below right] {$3$};
\draw (3,3) node[above left] {$1$};
\draw (5,2) node[below right] {$2$};
\draw (4,3) node[above left] {$2$};

\draw (6,2) node[below right] {$1$};
\draw (5,3) node[above left] {$3$};
\draw (6,3) node[below right] {$2$};
\draw (4,3) node[above left] {$2$};
\draw (5,4) node[above] {$1$};
\draw (7,3) node[right] {$1$};
\draw (6,4) node[above] {$1$};
\draw (7,4) node[above right] {$1$};

%-- (1,2) -- (3,2) -- (3,3) -- (5,3) -- (5,4) -- (6,4);

%\draw (7,4) node[above right] {$(7,4)$};
%\draw (0,0) node[below left] {$(0,0)$};
%\draw[fill=black](2,1) circle (0.5mm);

%\draw[fill=blue](1,4/7) circle (1mm); \draw[fill=red](7/4,1) circle (1mm); \draw[fill=blue](2,8/7) circle (1mm); 
%\draw[fill=blue](3,12/7) circle (1mm); \draw[fill=red](3.5,2) circle (1mm); \draw[fill=blue](4,16/7) circle (1mm);
%\draw[fill=blue](5,20/7) circle (1mm); \draw[fill=red](5.25,3) circle (1mm); \draw[fill=blue](6,24/7) circle (1mm);

\end{tikzpicture}

\caption{Paths of lower and upper Christoffel words $w=aabaabaabab$ and $\tilde w$ and the function $\delta_w$}
\label{Cword2}
\end{figure}

See for example Figure \ref{Cword2}.

\begin{theorem}\label{conjugate} Let $w=u^k$, $u$ primitive, $k\geq 2$. Then $w$ is factor-symmetric if and only if $u$ is the conjugate of some Christoffel word.
\end{theorem}

\noindent{\bf Open question}: which primitive trapezoidal words are strongly factor-symmetric? We know that if a word is primitive, Sturmian, and strongly factor-symmetric, it must be
a Christoffel word. Hence the question is really: which primitive trapezoidal words, that are not Sturmian, are strongly factor-symmetric? An example is the word $aabb$. The work of \cite{DA} might help.

\section{Sketch of proofs}

Proving Theorem \ref{ECwSS} amounts to proving that the bivariate (commutative) polynomial $\sum_{i,j}\delta_w(i,j)a^ib^j\in \mathbb N[a,b]$ is reciprocal, with an appropriate (but evident) 
definition of ``reciprocal". One shows that this property is preserved by product. Then one shows, using the notion of attractor and circular attractor \cite{MRRRS} 
that the factors of $w$, which intersect the cut in the palindromic factorization $uv$ of $w$, are all the factors of $w$; moreover they are distinct \cite{BR}. Hence 
the set of factors is the unambiguous product of the set of suffixes of $u$ by the set of prefixes of $v$. Making the letters commute, the previous polynomial is the 
product of two polynomials; these are reciprocal, by palindromicity of $u$ and $v$.

To prove Theorem \ref{SSSwC}, it is enough to prove that $w$ is unbordered. One show that $w$ has a nontrivial period $p$ if and only if the intersection of the 
support of $\delta_w$ and of the line of equation $x+y=p$ is a singleton. Hence, by symmetry of the support, if $w$ has this period, $w$ also has the period $n-p$, 
hence is not primitive, by a Fine-Wilf lemma.

Let us sketch the proof of Theorem \ref{SSiffuCw}. Suppose that $u$ is a Christoffel word. Clearly, all circular factors of $u$ are factors of $w$. An ad hoc construction then allows one to enumerate all factors of $w$, relating them to the circular factors of $u$, and implying that $\delta_w$ has the required symmetry.

Conversely, one shows that the hypothesis implies that $u$ has at most $k+1$ circular factors of length $k$, for $k=0,1,\ldots,|u|-1$; being primitive, it must have exactly $k+1$ factors. Hence $u$ is the conjugate of a Christoffel word. We conclude the result using periodicity as above. 

For Theorem \ref{bijection}, a closer look at the combinatorics behind the algebraic proof using polynomials gives the bijection. 
%One may even note that the bijection is decreasing with respect to the lexicographical order.

For Theorem \ref{support}, one notes that since $w$ is balanced, there at most two points in the intersection of the 
support of $\delta_w$ and the line of equation $x+y=p$. One of them is given by the intersection of the lower path and the line, and corresponds to the prefix of
length $p$ of $w$. The other to the suffix of length $p$ of $w$, since $w=amb$, $m$ palindrome.

Finally, if $w=u^k$ is factor symmetric, then one shows as above that $u$ is the conjugate of a Christoffel word. Conversely, each power of a conjugate of a Christoffel word is Sturmian, hence trapezoidal, hence factor-symmetric. This proves Theorem \ref{conjugate}.

%%%%%%%%%%%%%%%%%%%%%%%%%%%%%%%%%%%%%%%%%%%%%%%%%%%%%%%%%%%%%%%%%%%%%%%%%%%%%%%
%\nocite{*}
\bibliographystyle{eptcs}
\bibliography{main}

\end{document}